\documentclass[12pt]{article}
\usepackage{amssymb, amsfonts, amssymb, amsthm, amscd}

\newtheorem{theorem}{{\bf Theorem}}[section]
\newtheorem{lemma}{{\bf Lemma}}[section]
\newtheorem{proposition}{{\bf Proposition}}[section]

\newtheorem{remark}{{\bf Remark}}[section]

\begin{document}

\title{Reflected BSDE with quadratic growth and unbounded terminal value}
\author{J.-P. Lepeltier$^a$\thanks{
Email of corresponding author : lepeltier@univ-lemans.fr} and M. Xu$^{a,b}$%
\thanks{
Email : xvmingyu@gmail.com} \\
{\small {$^a$D\'epartement de Math\'ematiques, Universit\'e du Maine, Avenue
Olivier Messiaen,}}\\
{\small {\ 72085 Le Mans Cedex 9, France}}\\
{\small $^b$Department of Financial Mathematics and Control science, School
of Mathematical Science,} \\
{\small Fudan University, Shanghai, 200433, China.}}
\date{\small \textit{The first version of this paper is done in June, 2006, this version is in April, 2007}}\maketitle

\begin{center}
\textbf{Abstract }
\end{center}

{\small In this paper we prove the existence of a solution for reflected
BSDE's\ whose coefficient is of quadratic growth in $z$ and of linear growth
in $y$, with an unbounded terminal value.}

\vspace*{.8cm} \textbf{Keywords: }{\small Reflected Stochastic Differential
Equations, reflected ordinary differential equation, characterization of the
solution, quadratic growth}

\section{Introduction}

In this paper we are interested with the following real valued reflected
backward stochastic differential equations (RBSDE's in short) with one
continuous barrier

\begin{eqnarray*}
Y_{t}&=&\xi +\int_{t}^{T}f\left( s,Y_{s},Z_{s}\right)
ds-\int_{t}^{T}Z_{s}dB_{s}+K_{T}-K_{t},0\leq t\leq T \\
Y_{t}&\geq& L_{t},0\leq t\leq T,\int_{0}^{T}\left( Y_{s}-L_{s}\right)
dK_{s}=0
\end{eqnarray*}
where $\left( B_{t}\right) $ is a standard Brownian motion. In our setting
the coefficient, namely $f$, is of quadratic growth in $z$ and of linear
growth in $y$.

In 1996, El Karoui et al. \cite{EKPPQ} first introduced this kind of
equations and proved the existence and uniqueness of the solution under a
Lipschitz condition in $y$ and $z$. Then in 1997 Matoussi \cite{M97} studied
the case when $f$ is of linear growth in $y$ and $z$. When the terminal
value $\xi $ is square integrable he proved the existence of a maximal and a
minimal solution. Later RBSDE's, whose coefficients are quadratic growth in $%
z$, have been studied by Kobylanski, Lepeltier, Quenez, Torres in \cite{KLQT}%
, but they required the terminal value $\xi $ is bounded.

In an interesting paper, Briand and Hu \cite{BH2005} relaxed the boundness
of $\xi $ for non reflected BSDE's whose coefficients is quadratic growth in
$z$. In this work we use a similar approach in the case of RBSDE's, with the
help of existence results contained in \cite{KLQT}.

The next section is devoted to the assumptions and the claim of the main
result theorem 2.1. The third section gives some estimation results which
are important to establish the proof of theorem 2.1 in section 4. Then
section 5 is devoted to get an extension to the case that $f$ is superlinear
in $y$. Finally in section 6 (Appendix) we study the existence, uniqueness
and characterization of the solution for backward ordinary differential
equations with one lower continuous barrier, which is a key point in the
technics used in section 3.

\section{Assumptions and Main result}

Let $(\Omega ,\mathcal{F},P)$ be a complete probability space, and $%
(B_{t})_{0\leq t\leq T}=(B_{t}^{1},B_{t}^{2},\cdots ,B_{t}^{d})_{0\leq t\leq
T}^{\prime }$ be a $d$-dimensional Brownian motion defined on a finite
interval $[0,T]$, $0<T<+\infty $. Denote by $\{\mathcal{F}_{t};0\leq t\leq
T\}$ the standard filtration generated by the Brownian motion $B$, i.e. $%
\mathcal{F}_{t}$ is the completion of
\[
\mathcal{F}_{t}=\sigma \{B_{s};0\leq s\leq t\},
\]
with respect to $(\mathcal{F},P)$. We denote by $\mathcal{P}$ the $\sigma $%
-algebra of predictable sets on $[0,T]\times \Omega $.

We shall need the following spaces:

\[
\begin{array}{ll}
\mathbf{L}^{2}(\mathcal{F}_{t})= & \{\eta
:\mathcal{F}_{t}\mbox{-measurable
random real-valued variable, s.t. }E(|\eta |^{2})<+\infty \}, \\
\mathbf{H}_{n}^{2}(0,T)= & \{(\psi _{t})_{0\leq t\leq
T}:\mbox{predictable process valued in }\mathbb{R}^{n}\mbox{, s.t.
}E\int_{0}^{T}\left| \psi
(t)\right| ^{2}dt<+\infty \}, \\
\mathbf{S}^{2}(0,T)= & \{(\psi _{t})_{0\leq t\leq
T}:\mbox{progressively
measurable real-valued process,} \\
& \mbox{s.t. }E(\sup_{0\leq t\leq T}\left| \psi (t)\right|
^{2})<+\infty \},
\\
\mathbf{A}^{2}(0,T)= & \{(K_{t})_{0\leq t\leq T}:\ \mbox{adapted
continuous
increasing process, } \\
& \mbox{s.t. }K(0)=0\mbox{, }E(K(T)^{2})<+\infty \}.
\end{array}
\]
$\mathbf{S}^{\infty }\left( 0,T\right) $ denotes the set of predictable
bounded processes.

In this paper, we work under the following assumptions:

\noindent\textbf{Assumption 2.1.} a coefficient $f:[0,T]\times \mathbb{R\times R}%
^{d}\rightarrow \mathbb{R}$, is linear increasing in $y$ and quadratic
growth in $z$: there exists $\alpha $, $\beta \geq 0$, $\gamma >0$,
satisfying $\alpha \geq \frac{\beta }{\gamma }$, such that for $\forall
(t,y,z)\in [0,T]\times \mathbb{R}\times \mathbb{R}^{d}$,
\begin{equation}
\left| f(t,y,z)\right| \leq \alpha +\beta \left| y\right| +\frac{\gamma }{2}%
\left| z\right| ^{2};  \label{con-gro}
\end{equation}
moreover $f(t,y,z)$ is continuous in $(y,z)$, for all $t\in [0,T]$.

\noindent\textbf{Assumption 2.2.} a terminal condition $\xi \in \mathbf{L}^{2}(%
\mathcal{F}_{T})$, such that
\[
E[e^{\gamma e^{\beta T}\left| \xi \right| }]<+\infty .
\]

\noindent\textbf{Assumption 2.3.} a barrier $L$, which is a bounded
continuous process, with $L_{T}\leq \xi $, and for $\forall t\in
[0,T]$, $\left| L_{t}\right| \leq a_{t}$, where $a_{t}$ is a
deterministic and continuous process.

For the terminal condition, we propose another stronger assumption:

\noindent\textbf{Assumption 2.4 }a terminal time $\xi \in \mathbf{L}^{2}(\mathcal{F}%
_{T})$, such that $E[e^{2\gamma e^{\beta T}\left| \xi \right| }]<+\infty $.

\begin{remark}
\label{ass-ter}From the assumption 2.3, we know that $\xi $ has a lower
bound in view of $\xi \geq L_{T}\geq -a_{T}$.
\end{remark}

Our main result in this paper is:

\begin{theorem}
\label{main}Under the assumptions 2.1-2.3, the reflected BSDE associated to $%
(\xi ,f,L)$ admits at least a solution, i.e. there exists a triplet $%
(Y_{t},Z_{t},K_{t})_{0\leq t\leq T}$, with $Y\in \mathbf{S}^{2}(0,T)$, and $%
K\in \mathbf{A}^{2}(0,T)$, such that
\begin{eqnarray*}
Y_{t} &=&\xi
+\int_{t}^{T}f(s,Y_{s},Z_{s})ds+K_{T}-K_{t}-\int_{t}^{T}Z_{s}dB_{s}, \\
Y_{t} &\geq &L_{t},\;\;\int_{0}^{T}(Y_{t}-L_{t})dK_{t}=0.
\end{eqnarray*}
Moreover if assumption 2.4 holds, then $Z\in \mathbf{H}_{d}^{2}(0,T)$.
\end{theorem}

\section{Estimation results}

To prove theorem \ref{main}, we need prove an estimation result. Define
\[
\mathbf{L}_{\gamma }^{2}(\mathbb{R})=\{\mbox{ }(v_{t})_{0\leq t\leq
T}:[0,T]\rightarrow \mathbb{R}\mbox{, s.t. }\int_{0}^{T}e^{\gamma
t}\left| v_{t}\right| ^{2}dt<\infty \}\mbox{, for }\gamma \in
\mathbb{R}.
\]

\begin{lemma}
\label{est1}Let assumption 2.1 hold and $\xi $ be a bounded $\mathcal{F}_{T}$%
-measurable random variable. If $(Y_{t},Z_{t},K_{t})_{0\leq t\leq T}$ is a
solution of the RBSDE$(\xi ,f,a)$ in $\mathbf{S}^{\infty }(0,T)\times
\mathbf{H}_{d}^{2}(0,T)\times \mathbf{A}^{2}(0,T)$, then
\[
L_{t}\leq Y_{t}\leq \frac{1}{\gamma }\ln (E[\theta _{t}(\xi )|\mathcal{F}%
_{t}]).
\]
Here the mapping $\theta _{t}(\cdot ):\mathbf{R\rightarrow L}_{\gamma }^{2}(%
\mathbf{R)}$ is defined that for $x\in \mathbf{R}$, $(\theta
_{t}(x),k_{t}(x))$ is the unique solution of following reflected backward
ordinary differential equation,
\begin{eqnarray}
\theta _{t}(x) &=&e^{\gamma x}+\int_{t}^{T}H(\theta
_{s}(x))ds+k_{T}(x)-k_{t}(x),  \label{RBODE} \\
\theta _{t}(x) &\geq &e^{\gamma a_{t}},\int_{0}^{T}(\theta _{t}(x)-e^{\gamma
a_{t}})dk_{t}(x)=0.  \nonumber
\end{eqnarray}
with $H(p)=p(\alpha \gamma +\beta \ln p)1_{[1,+\infty )}(p)+\gamma \alpha
1_{(-\infty ,1)}(p)$.
\end{lemma}

\proof%
Consider the change of variable
\[
P_{t}=e^{\gamma Y_{t}},Q_{t}=\gamma e^{\gamma Y_{t}}Z_{t}=\gamma
P_{t}Z_{t},J_{t}=\int_{0}^{t}\gamma e^{\gamma Y_{s}}dK_{s}.
\]
It is easy to check that $(Y,Z,K)$ is a solution of the RBSDE$(\xi ,f,L)$ if
and only if $(P,Q,J)$ is a solution of the RBSDE$(e^{\gamma \xi
},F,e^{\gamma L_{t}})$, where
\[
F(s,p,q)=1_{\{p>0\}}(\gamma pf(s,\frac{\ln p}{\gamma },\frac{q}{\gamma p})-%
\frac{1}{2}\frac{\left| q\right| ^{2}}{p}),
\]
i.e. the triplet $(P_{t},Q_{t},J_{t})_{0\leq t\leq T}$ satisfies
\begin{eqnarray*}
P_{t} &=&e^{\gamma \xi
}+\int_{t}^{T}F(s,P_{s},Q_{s})ds+J_{T}-J_{t}-\int_{t}^{T}Q_{s}dB_{s}, \\
P_{t} &\geq &e^{\gamma L_{t}},\int_{0}^{T}(P_{t}-e^{\gamma L_{t}})dJ_{t}=0.
\end{eqnarray*}
Then in order to get the integral property of $Y$, it is sufficient to study
the integrability of the process $P$. First $P_{t}\geq e^{\gamma L_{t}}$,
then it remains to find out an upper bound of $P$.

We define the mapping $\theta _{t}(\cdot ):\mathbf{R\rightarrow
L}_{\gamma }^{2}(\mathbf{R)}$, for $x\in \mathbf{R}$, $\theta
_{t}(x)$ with an increasing process $k_{t}(x)$, is a unique
solution of the reflected BODE with coefficient $H$, deterministic
barrier $e^{\gamma a_{t}}$, and terminal condition $e^{\gamma
x}\in \mathbb{R}$, satisfying $x\geq a_{T}$; i.e. (\ref {RBODE})
is satisfied. Thanks to theorem \ref{ode-main} in the Appendix, we
know that $\theta _{t}(x)$ exists and can be written in the
following forms
\begin{eqnarray*}
\theta _{t}(x) &=&\sup_{t\leq s\leq T}\varphi _{t}(s,e^{\gamma
a_{s}}1_{\{s<T\}}+e^{\gamma x}1_{\{s=T\}}) \\
&=&\max \{\varphi _{t}(T,e^{\gamma x}),\sup_{t\leq s<T}\varphi
_{t}(s,e^{\gamma a_{s}})\} \\
&=&\sup_{t\leq s\leq T}[\int_{t}^{s}H(\theta _{r}(x))dr+e^{\gamma
a_{s}}1_{\{s<T\}}+e^{\gamma x}1_{\{s=T\}}],
\end{eqnarray*}
where $\varphi _{t}(s,e^{\gamma a_{s}}1_{\{s<T\}}+e^{\gamma x}1_{\{s=T\}})$
is the solution of the non-reflected BODE on $[0,s]$ with coefficient $H$
and terminal value $e^{\gamma a_{s}}1_{\{s<T\}}+e^{\gamma x}1_{\{s=T\}}$,
i.e. the followings hold
\begin{eqnarray}
\varphi _{t}(T,e^{\gamma x}) &=&e^{\gamma x}+\int_{t}^{T}H(\varphi
_{r}(T,e^{\gamma x}))dr,  \label{odes} \\
\varphi _{t}(s,e^{\gamma a_{s}}) &=&e^{\gamma
a_{s}}+\int_{t}^{s}H(\varphi _{r}(s,e^{\gamma a_{s}}))dr,\mbox{ for
}0\leq s<T.  \nonumber
\end{eqnarray}

For a bounded $\mathcal{F}_{T}$-measurable random variable $\xi $, we get
\[
\theta _{t}(\xi )=\max \{\varphi _{t}(T,e^{\gamma \xi }),\sup_{t\leq
s<T}\varphi _{t}(s,e^{\gamma a_{s}})\},
\]
which is also an $\mathcal{F}_{T}$-measurable random variable. Since
\[
\theta _{t}(\xi )=\sup_{t\leq s\leq T}[\int_{t}^{s}H(\theta _{r}(\xi
))dr+e^{\gamma a_{s}}1_{\{s<T\}}+e^{\gamma \xi }1_{\{s=T\}}],
\]
for any stopping time $\tau $, such that $t\leq \tau \leq T$, we have
\[
\theta _{t}(\xi )\geq \int_{t}^{\tau }H(\theta _{r}(\xi ))dr+e^{\gamma
a_{\tau }}1_{\{\tau <T\}}+e^{\gamma \xi }1_{\{\tau =T\}}.
\]
So
\[
\theta _{t}(\xi )\geq ess\sup_{\tau \in \mathcal{T}_{t,T}}\int_{t}^{\tau
}H(\theta _{r}(\xi ))dr+e^{\gamma a_{\tau }}1_{\{\tau <T\}}+e^{\gamma \xi
}1_{\{\tau =T\}},
\]
where $\mathcal{T}_{t,T}$ is the set of the stopping times valued in $[t,T]$.

Denote $\Theta _{t}(\xi ):=E[\theta _{t}(\xi )|\mathcal{F}_{t}]$, then we
have
\begin{eqnarray*}
\Theta _{t}(\xi ) &\geq &ess\sup_{\tau \in \mathcal{T}_{t,T}}E[\int_{t}^{%
\tau }H(\theta _{r}(\xi ))dr+e^{\gamma a_{\tau }}1_{\{\tau <T\}}+e^{\gamma
\xi }1_{\{\tau =T\}}|\mathcal{F}_{t}] \\
&\geq &ess\sup_{\tau \in \mathcal{T}_{t,T}}E[\int_{t}^{\tau }E[H(\theta
_{r}(\xi ))|\mathcal{F}_{r}]dr+e^{\gamma a_{\tau }}1_{\{\tau <T\}}+e^{\gamma
\xi }1_{\{\tau =T\}}|\mathcal{F}_{t}].
\end{eqnarray*}
Set $\Phi _{t}(\xi )$ equal to the right side; by the optimal stopping
problem, we know that there exist $(\Psi (\xi ),\Lambda (\xi ))\in \mathbf{H}%
_{d}^{2}(0,T)\times \mathbf{A}^{2}(0,T)$, such that $(\Phi (\xi ),\Psi (\xi
),\Lambda (\xi ))$ is the solution of the RBSDE$(e^{\gamma \xi },E[H(\theta
_{t}(\xi ))|\mathcal{F}_{t}],e^{\gamma a})$.

From assumption 2.1, it follows that the function $H$ is convex, increasing
in $p$. And $F(s,p,q)\leq H(p)$, for any $s\in [0,T]$, $p\in \mathbb{R}$, $%
q\in \mathbb{R}^{d}$. So for $r\in [0,T]$, we have
\begin{eqnarray*}
E[H(\theta _{r}(\xi ))|\mathcal{F}_{r}] &\geq &H(E[\theta _{r}(\xi )|%
\mathcal{F}_{r}])=H(\Theta _{r}(\xi )) \\
&\geq &H(\Phi _{r}(\xi ))\geq F(r,\Phi _{r}(\xi ),\Psi _{r}(\xi )).
\end{eqnarray*}

Since $\xi $ is a bounded $\mathcal{F}_{T}$-measurable random variable, it
follows that $\Phi _{t}(\xi )$ and $P_{t}$ are bounded. Since $H(p)$ is
locally Lipschitz, we can apply the trajectory comparison theorem for these
RBSDEs, and get for $t\in [0,T]$,
\[
\Theta _{t}(\xi )\geq \Phi _{t}(\xi )\geq P_{t}.
\]
Consequently
\[
Y_{t}\leq \frac{1}{\gamma }\ln \Theta _{t}(\xi )=\frac{1}{\gamma }\ln
(E[\theta _{t}(\xi )|\mathcal{F}_{t}]).
\]
$\square $

\begin{lemma}
\label{est2}Let assumptions 2.1 and 2.3 hold, and $\xi $ be a $\mathcal{F}%
_{T}$-measurable bounded random variable. If $(Y_{t},Z_{t},K_{t})_{0\leq
t\leq T}$ is a solution of the RBSDE$(\xi ,f,L)$, then
\begin{equation}
L_{t}\leq Y_{t}\leq \frac{1}{\gamma }\ln (E[\theta _{t}(\xi \vee a_{T})|%
\mathcal{F}_{t}]).  \label{est-bar}
\end{equation}
\end{lemma}

\proof%
Obviously $Y_{t}\geq L_{t}$. For the right side, consider the RBSDE$(\xi
\vee a_{T},f,a)$; since $a$ is a bounded continuous process, by \cite{KLQT},
it admits a maximal solution $(Y^{a},Z^{a},K^{a})$. From the comparison
theorem, we have $Y_{t}\leq Y_{t}^{a}$. Thanks to lemma \ref{est1}, $%
Y_{t}^{a}\leq \frac{1}{\gamma }\ln (E[\theta _{t}(\xi \vee a_{T})|\mathcal{F}%
_{t}])$, which follows
\[
Y_{t}\leq \frac{1}{\gamma }\ln (E[\theta _{t}(\xi \vee a_{T})|\mathcal{F}%
_{t}]).
\]
$\square $

\begin{remark}
\label{comp-s}We can also get some comparison results of $\theta
_{t}(x)$. Recalling the results in \cite{BH2005}, we can solve
equations (\ref{odes}) explicitly.  From their forms, it is easy to
check that $\varphi _{t}(T,e^{\gamma x})$ and $\varphi
_{t}(s,e^{\gamma a_{s}})$ are decreasing in $t$, and $\varphi
_{t}(T,e^{\gamma x})$ is increasing and continuous in $x $. So
$\theta _{t}(x)$ is increasing in $x$.

For $t_{1}$, $t_{2}\in [0,T]$, with $t_{1}\leq t_{2}$, we have
\[
\varphi _{t_{1}}(T,e^{\gamma x})\geq \varphi _{t_{2}}(T,e^{\gamma
x})\mbox{ and }\varphi _{t_{1}}(s,e^{\gamma a_{s}})\geq \varphi
_{t_{2}}(s,e^{\gamma a_{s}}).
\]
Remember that
\begin{eqnarray*}
\theta _{t_{1}}(x) &=&\max \{\varphi _{t_{1}}(T,e^{\gamma
x}),\sup_{t_{1}\leq s\leq t_{2}}\varphi _{t_{1}}(s,e^{\gamma
a_{s}}),\sup_{t_{2}\leq s\leq T}\varphi _{t_{1}}(s,e^{\gamma a_{s}})\}, \\
\theta _{t_{2}}(x) &=&\max \{\varphi _{t_{2}}(T,e^{\gamma
x}),\sup_{t_{2}\leq s\leq T}\varphi _{t_{2}}(s,e^{\gamma a_{s}})\},
\end{eqnarray*}
then we obtain $\theta _{t_{1}}(x)\geq \theta _{t_{2}}(x)$, i.e. $\theta
_{t}(x)$ is decreasing in $t$.
\end{remark}

\section{The proof of theorem \ref{main}}

Now we can prove our main result. Before beginning the proof, we present a
monotone stability theorem, which is proved in theorem 4 of \cite{KLQT}.

\begin{theorem}
\label{stab}Let $(\xi ^{p})_{p\in \mathbb{N}}$, $\xi $ be a family of
terminal condition, $(g^{p})_{p\in \mathbb{N}}$, $g$ be a family of
coefficients, $L$ be a continuous bounded process, which satisfy:

(a) there exists a constant $b>0$, such that for each $p$, $\left| \xi
^{p}\right| \leq b$, and $\left| L_{t}\right| \leq b$, for $t\in [0,T]$,
with $\xi ^{p}\geq L_{T}$.

(b) $g^{p}$, $g:[0,T]\times \Omega \times \mathbb{R}\times \mathbb{R}%
^{d}\rightarrow \mathbb{R}$ are $\mathcal{P}\otimes \mathcal{B}(\mathbb{R}%
)\otimes \mathcal{B}(\mathbb{R}^{d})$-measurable, and there exists a
function $l_{1}$ of the form $l_{1}(y)=a_{1}(1+\left| y\right| )$, with $%
a_{1}>0$, and a constant $A$, such that for each $p$,
\[
\left| g^{p}(t,y,z)\right| \leq l_{1}(y)+A\left| z\right| ^{2}\mbox{ and }%
\left| g(t,y,z)\right| \leq l_{1}(y)+A\left| z\right| ^{2}.
\]

(c) the sequence $(g^{p})$ converge increasingly (resp. decreasingly) to $g$
locally uniformly on $[0,T]\times \mathbb{R}\times \mathbb{R}^{d}$, and $%
(\xi ^{p})$ converge increasingly (resp. decreasingly) to $\xi $.

For each $p$, let $(Y^{p},Z^{p},K^{p})$ be the maximal solution of the RBSDE$%
(\xi ^{p},g^{p},L)$. Then the sequence $(Y^{p})$ converges increasingly
(resp. decreasingly) to $Y$ uniformly on $[0,T]$, $(Z^{p})$ converges to $Z$
in $\mathbf{H}_{d}^{2}(0,T)$, and $(K^{p})$ converges decreasingly (resp.
increasingly) to $K$ uniformly on $[0,T]$, where $(Y,Z,K)$ is the maximal
solution of the RBSDE$(\xi ,g,L)$.
\end{theorem}

\begin{remark}
The results still hold if we consider the minimal solutions of the RBSDEs.
\end{remark}

\textbf{Proof of theorem \ref{main}:}

By remark \ref{ass-ter}, we know that $\xi $ has a lower bound. So we only
need to consider the approximation of the upper side. For $n\geq a_{T}$, we
set $\xi ^{n}:=\xi \wedge n$. It is known from \cite{KLQT} that there exists
a maximal bounded solution $(Y^{n},Z^{n},K^{n})$ to the RBSDE$(\xi ^{n},f,L)$%
,
\begin{eqnarray*}
Y_{t}^{n} &=&\xi
^{n}+\int_{t}^{T}f(s,Y_{s}^{n},Z_{s}^{n})ds+K_{T}^{n}-K_{t}^{n}-%
\int_{t}^{T}Z_{s}^{n}dB_{s}, \\
Y_{t}^{n} &\geq &L_{t},\int_{0}^{T}(Y_{t}^{n}-L_{t})dK_{t}^{n}=0.
\end{eqnarray*}
Here $(Y^{n},Z^{n},K^{n})\in \mathbf{S}^{2}(0,T)\times \mathbf{H}%
_{d}^{2}(0,T)\times \mathbf{A}^{2}(0,T)$. Then from lemma \ref{est2}, we get
\[
L_{t}\leq Y_{t}^{n}\leq \frac{1}{\gamma }\ln (E[\theta _{t}(\xi ^{n}\vee
a_{T})|\mathcal{F}_{t}]).
\]
By the comparison theorem under superlinear condition in the Appendix of
\cite{X2004}, it follows that for $t\in [0,T]$,
\[
Y_{t}^{n}\leq Y_{t}^{n+1},K_{t}^{n}\geq K_{t}^{n+1}.
\]
Set $Y_{t}=\sup_{n}Y_{t}^{n}$, $K_{t}=\inf_{n}K_{t}^{n}$. By remark \ref
{comp-s}, we have $0\leq \theta _{t}(\xi ^{n}\vee a_{T})\leq \theta _{t}(\xi
\vee a_{T})\leq \theta _{0}(\xi \vee a_{T})$, then
\[
L_{t}\leq Y_{t}\leq \frac{1}{\gamma }\ln (E[\theta _{0}(\xi \vee a_{T})|%
\mathcal{F}_{t}])
\]
in view of the dominated convergence theorem and assumption 2.2. So $Y\in
\mathbf{S}^{2}(0,T)$ and $K\in \mathbf{A}^{2}(0,T)$ since $%
E[(K_{T})^{2}]\leq E[(K_{T}^{n})^{2}]$.

Let us introduce the following stopping times
\[
\tau _{m}=\inf \{t\in [0,T],\frac{1}{\gamma }\ln (E[\theta _{t}(\xi \vee
a_{T})|\mathcal{F}_{t}])\geq m\}\wedge T.
\]
Then denote $(Y^{n,m},Z^{n,m},K^{n,m})=(Y_{t\wedge \tau
_{m}}^{n},Z_{t}^{n}1_{\{t<\tau _{m}\}},K_{t\wedge \tau _{m}}^{n})$, which
satisfy the following RBSDE
\begin{eqnarray*}
Y_{t}^{n,m} &=&\xi ^{n,m}+\int_{t}^{T}1_{\{s\leq \tau
_{m}\}}f(s,Y_{s}^{n,m},Z_{s}^{n,m})ds+K_{T}^{n,m}-K_{t}^{n,m}-%
\int_{t}^{T}Z_{s}^{n,m}dB_{s}, \\
Y_{t}^{n,m} &\geq &L_{t},\int_{0}^{T}(Y_{t}^{n,m}-L_{t})dK_{t}^{n,m}=0,
\end{eqnarray*}
where $\xi ^{n,m}=Y_{T}^{n,m}=Y_{\tau _{m}}^{n}$.

For $m$ fixed, we have that $\{\xi ^{n,m}\}$ is increasing in $n$, and
bounded by $m$, in view of $\sup_{n}\sup_{t}\left| Y_{t}^{n,m}\right| \leq m$%
. Now we apply the monotone stability theorem \ref{stab} to $%
\{Y^{n,m}\}_{n\in \mathbb{N}}$. Setting $Y_{t}^{m}=\sup_{n}Y_{t}^{n,m}$,
then $Y^{n,m}$ converge uniformly to $Y^{m}$ on $[0,T]$ and there exist
processes $Z^{m}\in \mathbf{H}_{d}^{2}(0,T)$ and $K^{m}\in \mathbf{A}%
^{2}(0,T)$, such that $Z^{n,m}\rightarrow Z^{m}$ in $\mathbf{H}_{d}^{2}(0,T)$%
, and $K^{n,m}$ converges uniformly decreasingly to $K^{m}$. Furthermore, $%
(Y^{m},Z^{m},K^{m})$ solves
\begin{eqnarray*}
Y_{t}^{m} &=&\xi ^{m}+\int_{t}^{T}1_{\{s\leq \tau
_{m}\}}f(s,Y_{s}^{m},Z_{s}^{m})ds+K_{T}^{m}-K_{t}^{m}-%
\int_{t}^{T}Z_{s}^{m}dB_{s}, \\
Y_{t}^{m} &\geq &L_{t},\int_{0}^{T}(Y_{t}^{m}-L_{t})dK_{t}^{m}=0,
\end{eqnarray*}
where $\xi ^{m}=\sup_{n}Y_{\tau _{m}}^{n}$.

Since $\tau _{m}\leq \tau _{m+1}$, with the definition of $%
(Y^{m},Z^{m},K^{m})$, we deduce that
\[
Y_{t\wedge \tau _{m}}=Y_{t\wedge \tau
_{m}}^{m+1}=Y_{t}^{m},Z_{t}^{m+1}1_{\{t\leq \tau
_{m}\}}=Z_{t}^{m},K_{t\wedge \tau _{m}}=K_{t\wedge \tau
_{m}}^{m+1}=K_{t}^{m}.
\]
Since $Y^{m}$ and $K^{m}$ are continuous, and $P-a.s.$ $\tau _{m}=T$ for $m$
large enough , so $Y$ and $K$ are continuous on $[0,T]$. We define $Z$ on $%
[0,T)$ by setting
\[
Z_{t}=Z_{t}^{1}1_{\{t\leq \tau _{1}\}}+\sum_{m\geq 2}Z_{t}^{m}1_{(\tau
_{m-1},\tau _{m}]}(t),
\]
so $Z_{t}1_{\{t\leq \tau _{m}\}}=Z_{t}^{m}1_{\{t\leq \tau _{m}\}}=Z_{t}^{m}$%
and the triplet $(Y,Z,K)$ satisfies
\begin{equation}
Y_{t\wedge \tau _{m}}=Y_{\tau _{m}}+\int_{t\wedge \tau _{m}}^{\tau
_{m}}f(s,Y_{s},Z_{s})ds+K_{\tau _{m}}-K_{t\wedge \tau _{m}}-\int_{t\wedge
\tau _{m}}^{\tau _{m}}Z_{s}dB_{s}.  \label{equ-m}
\end{equation}
Since
\begin{eqnarray*}
P(\int_{0}^{T}\left| Z_{s}\right| ^{2}ds) &=&P(\int_{0}^{T}\left|
Z_{s}\right| ^{2}ds=\infty ,\tau _{m}=T)+P(\int_{0}^{T}\left| Z_{s}\right|
^{2}ds=\infty ,\tau _{m}<T) \\
&\leq &P(\int_{0}^{T}\left| Z_{s}\right| ^{2}ds=\infty )+P(\tau _{m}<T),
\end{eqnarray*}
with $\tau _{m}\nearrow T$, as $m\rightarrow \infty $, we deduce that $%
\int_{0}^{T}\left| Z_{s}\right| ^{2}ds<\infty $, $P$-a.s. Finally letting $%
m\rightarrow \infty $ in (\ref{equ-m}), we get that $(Y,Z,K)$ verifies the
equation.

On the other hand, $Y_{t}^{m}\geq L_{t}$, so $Y_{t}\geq L_{t}$ on $[0,T]$
and for each $m$, $\int_{0}^{T}(Y_{t}^{m}-L_{t})dK_{t}^{m}=0$, which implies
$\int_{0}^{\tau _{m}}(Y_{t}-L_{t})dK_{t}=0$, for each $m$. Furthermore $P$%
-a.s. for $m$ large enough, $\tau _{m}=T$ so. we have $%
\int_{0}^{T}(Y_{t}-L_{t})dK_{t}=0$, $P$-a.s..

To complete the proof, we need to prove that under the assumption 2.4 the
process $Z$ is in $\mathbf{H}_{d}^{2}(0,T)$.

If $(Y,Z,K)$ is a solution of the RBSDE$(\xi ,f,L)$ constructed as before,
then
\begin{equation}
L_{t}\leq Y_{t}\leq \frac{1}{\gamma }\ln (E[\theta _{t}(\left| \xi \right|
\vee a_{T})|\mathcal{F}_{t}]),E[(K_{T})^{2}]<+\infty .  \label{est-yk}
\end{equation}
So under the assumption 2.4, we get,
\begin{equation}
E[\sup_{0\leq t\leq T}e^{2\gamma \left| Y_{t}\right| }]<+\infty .
\label{est-p}
\end{equation}
For $n\geq 1$, let $\sigma _{n}$ be the following stopping time:
\[
\sigma _{n}=\inf \{t\geq 0,\int_{0}^{t}e^{2\gamma \left| Y_{s}\right|
}\left| Z_{s}\right| ^{2}\geq n\}\wedge T,
\]
and consider the following function
\[
v(x)=\frac{1}{\gamma ^{2}}(e^{\gamma x}-1-\gamma x).
\]

By It\^{o}'s formula applied to $v(\left| Y_{t}\right| )$, with the notation
\[
sgn(x)=\left\{
\begin{array}{c}
1,x>0, \\
-1,x\leq 0,
\end{array}
\right.
\]
we get on $[0,t\wedge \sigma _{n}]$,
\begin{eqnarray*}
v(\left| Y_{0}\right| ) &=&v(\left| Y_{t\wedge \sigma _{n}}\right|
)+\int_{0}^{t\wedge \sigma _{n}}[v^{\prime }(\left| Y_{s}\right|
)sgn(Y_{s})f(s,Y_{s},Z_{s})-\frac{1}{2}v^{\prime \prime }(\left|
Y_{s}\right| )\left| Z_{s}\right| ]ds \\
&&+\int_{0}^{t\wedge \sigma _{n}}v^{\prime }(\left| Y_{s}\right|
)sgn(Y_{s})dK_{s}-\int_{0}^{t\wedge \sigma _{n}}v^{\prime }(\left|
Y_{s}\right| )sgn(Y_{s})Z_{s}dB_{s}.
\end{eqnarray*}
From the assumption 2.1 and $v^{\prime }(x)\geq 0$, for $x>0$, we get
\begin{eqnarray}
v(\left| Y_{0}\right| ) &\leq &v(\left| Y_{t\wedge \sigma _{n}}\right|
)+\int_{0}^{t\wedge \sigma _{n}}v^{\prime }(\left| Y_{s}\right| )(\alpha
+\beta \left| Y_{s}\right| )ds+\sup_{0\leq s\leq T}(v^{\prime }(\left|
Y_{s}\right| )\cdot K_{T}  \label{inequ} \\
&&-\int_{0}^{t\wedge \sigma _{n}}v^{\prime }(\left| Y_{s}\right|
)sgn(Y_{s})Z_{s}dB_{s}-\frac{1}{2}\int_{0}^{t\wedge \sigma _{n}}(v^{\prime
\prime }(Y_{s})-\gamma v^{\prime }(\left| Y_{s}\right| ))\left| Z_{s}\right|
^{2}ds.  \nonumber
\end{eqnarray}
Notice that $(v^{\prime \prime }-\gamma v^{\prime })(x)=1$, for $x\geq 0$;
taking expectation in (\ref{inequ}), we get
\begin{eqnarray}
\frac{1}{2}E\int_{0}^{t\wedge \sigma _{n}}\left| Z_{s}\right| ^{2}ds &\leq
&E[\frac{1}{\gamma ^{2}}\sup_{0\leq s\leq T}e^{\gamma \left| Y_{t}\right| }+%
\frac{1}{\gamma }\int_{0}^{T}e^{\gamma \left| Y_{t}\right| }(\alpha +\beta
\left| Y_{s}\right| )ds]  \label{est-z} \\
&&+\frac{1}{\gamma }(E[\sup_{0\leq s\leq T}e^{2\gamma \left| Y_{s}\right|
}])^{\frac{1}{2}}\cdot (E[(K_{T})^{2}])^{\frac{1}{2}}.  \nonumber
\end{eqnarray}
By Fatou's lemma, with (\ref{est-yk}) and (\ref{est-p}), letting $%
n\rightarrow \infty $ in (\ref{est-z}), we obtain $E\int_{0}^{T}\left|
Z_{s}\right| ^{2}ds<\infty $. $\square $

\section{One extension}

In this section, we extend our results to a more general case when the
coefficient $f$ is superlinear in $y$. Let $h:\mathbb{R}_{+}\rightarrow %
\mathbb{R}_{+}$ be a non-decreasing convex $\mathcal{C}^{1}$ function with $%
h(0)>0$ such that
\begin{equation}
\int_{0}^{+\infty }\frac{du}{h(u)}=+\infty ,\mbox{ and
}\sup_{y>0}e^{-\gamma y}h(y)<+\infty  \label{slin}
\end{equation}
We assume:

\noindent\textbf{Assumption 2.5.} the coefficient $f$ is continuous in $(y,z)$ for $%
t\in [0,T]$, and there exists $\gamma >0$ such that for $(t,y,z)\in
[0,T]\times \mathbb{R\times R}^{d}$,
\[
\left| f(t,y,z)\right| \leq h(\left| y\right| )+\frac{\gamma }{2}\left|
z\right| ^{2}.
\]
Obviously, the linear increasing condition in assumption 2.1 corresponds to $%
h(y)=\alpha +\beta y$, but we can also give a superlinear growth in $y$, for
example we can take $h(y)=\alpha (y+e)\ln (y+e)$.

Before giving our integrability condition for the terminal value $\xi $, we
need some modifications. According to (\ref{slin}), we denote $%
c_{0}=\sup_{p\in (0,1)}\gamma ph(-\frac{\ln p}{\gamma })$ and
\[
p_{0}=\inf \{p\geq 1:\gamma ph(\frac{\ln p}{\gamma })\geq c_{0}\}.
\]
Finally, we define
\[
H(p)=\gamma ph(\frac{\ln p}{\gamma })1_{\{p\geq
p_{0}\}}+c_{0}1_{\{p<p_{0}\}}.
\]
Then $H$ is convex and we have the following lemma.

\begin{lemma}
For $x\in \mathbb{R}$, the reflected BODE
\begin{eqnarray*}
\theta _{t}(x) &=&e^{\gamma x}+\int_{t}^{T}H(\theta
_{s}(x))ds+k_{T}(x)-k_{t}(x), \\
\theta _{t}(x) &\geq &a_{t},\int_{0}^{T}(\theta _{t}(x)-a_{t})dk_{t}(x)=0.
\end{eqnarray*}
has a unique solution $(\theta _{t}(x),k_{t}(x))_{0\leq t\leq T}$. Moreover $%
\theta _{t}(x)$ is decreasing on $t$ and continuous increasing on $x$.
\end{lemma}

\proof%
The results follows easily from the representation of the solution:
\begin{eqnarray*}
\theta _{t}(x) &=&\max \{\varphi _{t}(T,e^{\gamma x}),\sup_{t\leq
s<T}\varphi _{t}(s,e^{\gamma a_{s}})\} \\
&=&\sup_{t\leq s\leq T}[\int_{t}^{s}H(\theta _{r}(x))dr+e^{\gamma
a_{s}}1_{\{s<T\}}+e^{\gamma x}1_{\{s=T\}}],
\end{eqnarray*}
where $\varphi _{t}(s,e^{\gamma a_{s}})$ (resp. $\varphi _{t}(T,e^{\gamma
x}) $) is a solution of ODE on $[0,s]$ (resp. $[0,T]$) associated to $%
(e^{\gamma a_{s}},H)$ (resp. $(e^{\gamma x},H)$), and the existence results
about the non reflected ODE$(e^{\gamma a_{s}},H)$, see lemma 6 in \cite
{BH2005}. $\square $

Now we give our third integrability condition for the terminal condition $%
\xi $:

\noindent\textbf{Assumption 2.6. }$\theta _{0}(\xi \vee a_{T})$ is
integrable.

Exactly as in the linear case, we can prove the following existence result:

\begin{theorem}
Under assumptions 2.3, 2.5 and 2.6, the reflected BSDE associated to $(\xi
,f,L)$ has at least one solution $(Y,Z,K)$ such that
\begin{eqnarray*}
Y_{t} &=&\xi
+\int_{t}^{T}f(s,Y_{s},Z_{s})ds+K_{T}-K_{t}-\int_{t}^{T}Z_{s}dB_{s}, \\
Y_{t} &\geq &L_{t},\;\;\int_{0}^{T}(Y_{t}-L_{t})dK_{t}=0.
\end{eqnarray*}
Moreover, we have $L_{t}\leq Y_{t}\leq \frac{1}{\gamma }\ln (E[\theta
_{t}(\xi \vee a_{T})|\mathcal{F}_{t}]).$
\end{theorem}

\section{Appendix}

\subsection{Trajectory comparison theorem}

In this subsection, we prove a trajectory comparison theorem for RBSDE's
under a Lipschitz condition.

\begin{theorem}
Suppose that for $i=1,2$, $\xi ^{i}\in \mathbf{L}^{2}(\mathcal{F}_{T})$, $%
f^{i}(t,y,z)$ are Lipschitz functions in $y$ and $z$, i.e. there exists a $%
\mu >0$, such that for $y_{1}$, $y_{2}\in \mathbb{R}$, $z_{1}$, $z_{2}\in %
\mathbb{R}^{d}$,
\[
\left| f^{i}(t,y_{1},z_{1})-f^{i}(t,y_{2},z_{2})\right| \leq \mu (\left|
y_{1}-y_{2}\right| +\left| z_{1}-z_{2}\right| ),
\]
with $f^{i}(t,0,0)\in \mathbf{H}^{2}(0,T)$, and $L^{i}$ are adapted
continuous processes, with
\[
\xi ^{i}\geq L_{T}^{i}\mbox{ and
}E(\sup_{t}((L_{t}^{i})^{+})^{2})<+\infty .
\]
Let $(Y^{i},Z^{i},K^{i})$, $i=1,2$, be the solutions of the RBSDE's$(\xi
^{i},f^{i},L^{i})$, respectively. Moreover, we set $\forall t\in [0,T]$, $P$%
-a.s.
\[
\xi ^{1}\leq \xi ^{2},f^{1}(t,Y_{t}^{1},Z_{t}^{1})\leq
f^{2}(t,Y_{t}^{1},Z_{t}^{1}),L_{t}^{1}\leq L_{t}^{2}.
\]
Then $Y_{t}^{1}\leq Y_{t}^{2}$, a.s., for $t\in [0,T]$.
\end{theorem}

\begin{remark}
We have the same result under the condition $f^{1}(t,Y_{t}^{2},Z_{t}^{2})%
\leq f^{2}(t,Y_{t}^{2},Z_{t}^{2})$.
\end{remark}

\proof%
Applying It\^{o}'s formula to $[(Y_{t}^{1}-Y_{t}^{2})^{+}]^{2}$,
then taking expectation, with Lipschitz condition, we get
\[
E[(Y_{t}^{1}-Y_{t}^{2})^{+}]^{2}\leq (2\mu ^{2}+2\mu
)E\int_{t}^{T}[(Y_{s}^{1}-Y_{s}^{2})^{+}]^{2}ds.
\]
From Gronwall's inequality, we deduce that $(Y_{t}^{1}-Y_{t}^{2})^{+}=0$, $%
t\in [0,T]$, i.e. $Y_{t}^{1}\leq Y_{t}^{2}$. $\square $

\subsection{Existence and uniqueness of a solution for reflected backward
ODE's with one continuous barrier}

We recall the definition of the space
\[
\mathbf{L}_{\gamma }^{2}(\mathbb{R})=\{\mbox{ }(v_{t})_{0\leq t\leq
T}:[0,T]\rightarrow \mathbb{R}\mbox{, s.t. }\int_{0}^{T}e^{\gamma
t}\left| v_{t}\right| ^{2}dt<\infty \}\mbox{, for }\gamma \in
\mathbb{R}.
\]
Consider the reflected backward ordinary differential equation(reflected
BODE in short) reflected to one continuous barrier $l$ on $[0,T]$, with
terminal value $x\in \mathbb{R}$, whose solution is a couple $%
(y_{t},k_{t})_{0\leq t\leq T}$, with $y\in \mathbf{L}_{0}^{2}(\mathbb{R})$
is continuous, and $k$ is a continuous increasing process, $k_{0}=0$, and
the followings hold
\begin{eqnarray}
y_{t} &=&x+\int_{t}^{T}\phi (y_{s})ds+k_{T}-k_{t},  \label{RODE} \\
y_{t} &\geq &l_{t},\int_{0}^{T}(y_{s}-l_{s})dk_{s}=0.  \nonumber
\end{eqnarray}
Here we suppose

\noindent\textbf{Assumption A1. }the function $\phi
:\mathbb{R\rightarrow R}$, is continuous, and there exists a
strictly positive function $l_{0}$, such that
$\left| \phi (y)\right| \leq l_{0}(y)$, with $\int_{0}^{\infty }\frac{dy}{%
l_{0}(y)}=\infty $. And $\phi $ is increasing in $y$.

\noindent\textbf{Assumption A2. }the barrier $l$ satisfies: $\alpha
\leq l_{t}\leq \beta $, with $\beta >1$, $0<\alpha \leq 1$. And
$l_{T}\leq x$.

Furthermore we assume that

\noindent\textbf{Assumption A3. }the non reflected BODE's with any terminal value $x$%
, any terminal time $0\leq s\leq T,$ and the coefficient $\phi $, have a
unique solution.

Our main result is the

\begin{theorem}
\label{ode-main}Under assumptions A1, A2 and A3, the reflected BODE (\ref
{RODE}) admits one unique solution $(y_{t},k_{t})_{0\leq t\leq T}$. Moreover
\[
y_{t}=\sup_{t\leq s\leq T}u_{t}^{s}=\sup_{t\leq s\leq T}[\int_{t}^{s}\phi
(y_{r})dr+x1_{\{s=T\}}+l_{s}1_{\{s<T\}}],
\]
where $(u_{t}^{s})_{0\leq t\leq s}$ is the unique solution of the following
ODE defined on $[0,s]$%
\[
u_{t}^{s}=(x1_{\{s=T\}}+l_{s}1_{\{s<T\}})+\int_{t}^{s}\phi (u_{r}^{s})dr.
\]
\end{theorem}

\begin{remark}
The solution $y$ is the smallest process which saisfies the equation and $%
y_{t}\geq l_{t}$, $t\in [0,T]$, i.e. if another couple $(y^{\prime
},k^{\prime })$ satisfies aussi the equation and $y_{t}^{\prime }\geq l_{t}$%
, then $y_{t}\leq y_{t}^{\prime }$. But the increasing process $k$ is not
the smallest one.
\end{remark}

We first consider the existence of a solution.

\subsubsection{Existence}

For the existence, we do not need the monotonicity condition of $\phi $ in $%
y $ in assumption A1 and assumption A3. The proof is done in three steps:

a)$\phi $ is Lipschitz in $y$,

b) $\phi $ is linear increasing in $y$,

c) $\phi $ is superlinear increasing in $y$.

We consider first

\textbf{a)} the case $\phi $ Lipschitz in $y$, i.e. there exists a constant $%
\mu \in \mathbb{R}$, such that for $y$, $y^{\prime }\in \mathbb{R}$, $\left|
\phi (y)-\phi (y^{\prime })\right| \leq \mu \left| y-y^{\prime }\right| $.

When $\phi =\phi _{t}$ in $\mathbf{L}_{0}^{2}(\mathbb{R)}$, which means $%
\phi $ does not depend of $y$, it is easy to check that the solution of such
an equation is $y_{t}=\max \{x+\int_{t}^{T}\phi _{s}ds,l_{t}\}$, $%
k_{t}=\int_{0}^{t}(l_{s}-(x+\int_{s}^{T}\phi _{r}dr))^{+}ds$. Thanks
to the Lipschitz property of $\phi $ we can construct a strict contraction in $\mathbf{L}_{\gamma }^{2}(\mathbb{R}%
{)}$, beginning with a given process $\{y^{1}\}\in
\mathbf{L}_{0}^{2}(\mathbb{R)}$. So the reflected BODE admits one
unique solution. $\square$

Moreover, we have a comparison theorem:
\begin{theorem}
\label{compODE1}We consider the equations associated to $(x^{i},\phi ^{i},l)$%
, $i=1,2$, and assume that $\phi ^{1}$ and $\phi ^{2}$ satisfy the Lipschitz
assumptions. Let $(y^{i},k^{i})$ be the respective solutions of these
equations. Moreover, we assume for $t\in [0,T]$,
\[
x^{1}\geq x^{2},\phi ^{1}(y_{t}^{1})\geq \phi ^{2}(y_{t}^{1}),l_{t}^{1}\geq
l_{t}^{2}.
\]
Then $y_{t}^{1}\geq y_{t}^{2}$.
\end{theorem}

\textbf{Proof.} We consider $((y_{t}^{2}-y_{t}^{1})^{+})^{2}$.
Notice that on the set $\{y_{t}^{2}\geq y_{t}^{1}\}$, $y_{t}^{2}>y_{t}^{1}%
\geq l_{t}^{1}\geq l_{t}^{2}$, so we have
\begin{eqnarray*}
&&\int_{t}^{T}(y_{t}^{2}-y_{t}^{1})^{+}d(k_{s}^{2}-k_{s}^{1}) \\
&\leq
&\int_{t}^{T}(y_{t}^{2}-l_{t}^{2})dk_{s}^{2}-%
\int_{t}^{T}(y_{t}^{1}-l_{t}^{2})1_{\{y_{t}^{2}>y_{t}^{1}\}}dk_{s}^{2}-%
\int_{t}^{T}(y_{t}^{2}-y_{t}^{1})^{+}dk_{s}^{1}\leq 0.
\end{eqnarray*}
Consequently, we get
\[
((y_{t}^{2}-y_{t}^{1})^{+})^{2}\leq
2k\int_{t}^{T}((y_{s}^{2}-y_{s}^{1})^{+})^{2}ds.
\]
It follows immediately that $(y_{t}^{2}-y_{t}^{1})^{+}=0$, i.e. $%
y_{t}^{1}\geq y_{t}^{2}$. $\square $

\begin{remark}
The result is still true under the assumption $\phi ^{1}(y_{t}^{2})\geq \phi
^{2}(y_{t}^{2})$, $t\in [0,T]$.
\end{remark}

\textbf{b)} We now suppose that $\phi $ is continuous and linear increasing
in $y$, i.e. there exists a constant $\mu _{l}\in \mathbb{R}$, such that for
$y\in \mathbb{R}$, $\left| \phi (t,y)\right| \leq \mu _{l}(1+\left| y\right|
)$.

\begin{lemma}
Under the assumptions \textbf{b)} and A2, there exists a minimal solution $%
(y_{t},k_{t})_{0\leq t\leq T}$ of the reflected BODE$(x,\phi ,l)$.
\end{lemma}

\proof%
We consider the following approximation: for $n\in \mathbb{N}$, define
\begin{equation}
\phi _{n}(y)=\inf_{x\in \mathbf{Q}}\{\phi (x)+n\left| y-x\right| \},
\label{min-app}
\end{equation}
then for $n\geq \mu _{l}$, $\phi _{n}$ satisfies
\begin{equation}
\begin{array}{ll}
1) & \mbox{Linear increasing: }\left| \phi _{n}(y)\right| \leq \mu
_{l}(1+\left| y\right| ); \\
2) & \mbox{Monotonicity: }\phi _{n}(y)\nearrow \phi (y); \\
3) & \mbox{Lipschitz condition: }\left| \phi _{n}(y)-\phi
_{n}(y^{\prime
})\right| \leq n\left| y-y^{\prime }\right| ; \\
4) & \mbox{Strong convergence: If }y_{n}\rightarrow y,\mbox{ then
}\phi _{n}(y_{n})\rightarrow \phi (y)\mbox{, as }n\rightarrow \infty
.
\end{array}
\label{con-phi}
\end{equation}

By the result of a), for each $n\in \mathbb{N}$, there exists a unique
solution $(y^{n},k^{n})$ of the equation $(x,\phi _{n},l)$. It's easy to
check that the solutions $(y^{n})$ are bounded uniformly in $n$, i.e. $%
\sup_{0\leq t\leq T}(y_{t}^{n})^{2}\leq C$. Thanks to the comparison theorem
\ref{compODE1}, and 2) of (\ref{con-phi}), we know that $y_{t}^{n}\nearrow
y_{t}$, for $t\in [0,T]$. By Fatou's lemma, we get $\sup_{0\leq t\leq
T}(y_{t})^{2}\leq C$, and $\int_{0}^{T}\left| y_{t}^{n}-y_{t}\right|
^{2}ds\rightarrow 0$, in view of the dominated convergence theorem.

Then we prove that the convergence still holds in some stronger sense; for $%
n $, $p\in \mathbb{N}$, we have
\[
\sup_{0\leq t\leq T}(y_{t}^{n}-y_{t}^{p})^{2}\leq
2(\int_{0}^{T}(y_{s}^{n}-y_{s}^{p})^{2}ds)^{\frac{1}{2}}(\int_{0}^{T}(\phi
_{n}(y_{s}^{n})-\phi _{p}(y_{s}^{p}))^{2}ds)^{\frac{1}{2}}.
\]
By 1) of (\ref{con-phi}) and the estimate of $(y^{n})$, we get easily $%
\int_{0}^{T}(\phi _{n}(y_{s}^{n})-\phi _{p}(y_{s}^{p}))^{2}ds\leq C$, so $%
\sup_{0\leq t\leq T}(y_{t}^{n}-y_{t}^{p})^{2}\rightarrow 0$, as $n$, $%
p\rightarrow \infty $ and the limit $y$ is continuous.

For $\{k^{n}\}$, it is easy to check that $\sup_{0\leq t\leq
T}(k_{t}^{n}-k_{t}^{p})^{2}\rightarrow 0$, as $n$, $p\rightarrow \infty $.
Then there exists a increasing continuous process $k$, such that $(y,k)$
satisfies the equation. At last, we consider
\[
\int_{0}^{T}(y_{t}^{n}-l_{t})dk_{t}^{n}-\int_{0}^{T}(y_{t}-l_{t})dk_{t}
\leq \int_{0}^{T}(y_{t}-l_{t})d(k_{t}^{n}-k_{t})\rightarrow 0,
\]
in view of $\sup_{0\leq t\leq T}(k_{t}^{n}-k_{t})^{2}\rightarrow 0$, as $%
n\rightarrow \infty $. Since $y_{t}^{n}\geq l_{t}$, we get $y_{t}\geq l_{t}$%
, for $t\in [0,T]$. With $\int_{0}^{T}(y_{t}^{n}-l_{t})dk_{t}^{n}=0$, we
have $\int_{0}^{T}(y_{t}-l_{t})dk_{t}=0$. The proof is complete. $\square $

For the maximal solution, it is sufficient to replace (\ref{min-app}) by
\[
\phi _{n}(y)=\sup_{x\in \mathbf{Q}}\{\phi (x)-n\left| y-x\right| \},
\]
which is a sequence of Lipschitz functions which converge decreasingly to $%
\phi $. Then using the same approximation method as before, we obtain the
existence of the maximal solution. We have also the following comparison
theorem.

\begin{theorem}
Let us \label{compODE2}consider $\phi _{1}$, $\phi _{2}$ which satisfy the
condition \textbf{b)}. We suppose that for $y\in \mathbb{R}$, $t\in [0,T],$%
\[
x^{1}\geq x^{2},\phi _{1}(y)\geq \phi _{2}(y),l_{t}^{1}\geq l_{t}^{2}.
\]
For the maximal (minimal) solution $(y^{i},k^{i})$, $i=1,2$, of the
reflected equation associated to $(x^{i},\phi ^{i},l)$, we have $%
y_{t}^{1}\geq y_{t}^{2}$, for $t\in [0,T]$.
\end{theorem}

\proof%
The result comes easily from the approximation and theorem \ref{compODE1}. $%
\square $

\textbf{c)}We consider finally the case $\phi $ is continuous and
superlinear in $y$, i.e. $\left| \phi (y)\right| \leq l_{0}(y)$, with $%
\int_{0}^{\infty }\frac{dy}{l_{0}(y)}=\infty $.

Let $v_{t}$ be the solution of the ordinary differential equation: $%
v_{t}=b+\int_{t}^{T}l_{0}(v_{s})ds$, where $b=x\vee \sup_{t}l_{t}$ (see \cite
{KLQT} Lemma 1). Then we have

\begin{lemma}
Under the assumption \textbf{c)}, the reflected equation has a maximal
solution $(y,k)$, which satisfies: $\underline{m}\leq l_{t}\leq y_{t}\leq
v_{t}\leq v_{0}$ where $\underline{m}:=\inf_{t}l_{t}$.
\end{lemma}

\proof%
Let $\rho :\mathbb{R}_{+}$ $\rightarrow \mathbb{R}_{+}$ be a smooth function
such that
\[
\rho (x)=\left\{
\begin{array}{ll}
\frac{r}{2}, & 0<x<\frac{r}{2}; \\
x, & r\leq x\leq R; \\
2R, & x>2R.
\end{array}
\right.
\]
Here $r$ and $R$ are two real number such that $0<r<\underline{m}$ and $%
R>v_{0}$. It is a direct result that the unique solution of the equation
\[
v_{t}^{\rho }=b+\int_{t}^{T}l_{0}(\rho (v_{s}^{\rho }))ds
\]
satisfies $v_{t}^{\rho }=v_{t}$ and $v_{t}^{\rho }=v_{t}\geq v_{T}\geq l_{t}$
for $t\in [0,T]$. So $(v^{\rho },0)$ can be considered as the solution of
the reflected BODE associated to $(b,l_{0}(y),l)$. Then we consider the
following reflected BODE with one barrier $l$:
\begin{eqnarray*}
y_{t}^{\rho } &=&x+\int_{t}^{T}\phi (\rho (y_{s}^{\rho }))ds+k_{T}^{\rho
}-k_{t}^{\rho }, \\
y_{t}^{\rho } &\geq &l_{t},\int_{0}^{T}(y_{t}^{\rho }-l_{t})dk_{t}^{\rho }=0.
\end{eqnarray*}
Since $\phi (\rho (y))$ is bounded and continuous, this equation admits a
maximal solution $(y^{\rho },k^{\rho })$. Thanks to the comparison theorem,
we get $y_{t}^{\rho }\leq v_{t}^{\rho }\leq v_{0}<R$. With $y_{t}^{\rho
}\geq l_{t}\geq \underline{m}>r$, it follows that
\[
\phi (\rho (y_{s}^{\rho }))=\phi (y_{s}^{\rho }),
\]
i.e. $(y^{\rho },k^{\rho })$ is also a maximal solution of the reflected
BODE associated to $(x,\phi ,l)$. $\square $

We have still a comparison theorem, which follows easily from the proof of
existence and theorem \ref{compODE2}.

\begin{theorem}
\label{compODE3}Consider functions $\phi _{1}$, $\phi _{2}$ which satisfy
condition \textbf{c)}. We suppose, for $(t,y)\in [0,T]\times \mathbb{R}$.
\[
x^{1}\geq x^{2},\phi _{1}(y)\geq \phi _{2}(y),l_{t}^{1}\geq l_{t}^{2},
\]
For the maximal (minimal) solutions $(y^{i},k^{i})$, $i=1,2$, of the
reflected equations associated to $(x^{i},\phi ^{i},l)$, we have $%
y_{t}^{1}\geq y_{t}^{2}$, for $t\in [0,T]$.
\end{theorem}

\subsubsection{Uniqueness and characterization of the solution}

Here we will give a characterization of the solution of the reflected BODE
under assumptions A1, A2 and A3. First, we consider the following lemma:

\begin{lemma}
\label{ode-lim}Let $u^{\varepsilon }$ be the unique solution of the
following BODE for some $\varepsilon >0$, $\varepsilon \in \mathbb{R}$,
\begin{equation}
u_{t}^{\varepsilon }=x-\varepsilon +\int_{t}^{T}\phi (u_{s}^{\varepsilon
})ds,  \label{ode-epsilon}
\end{equation}
Then $u_{t}^{\varepsilon }$ converge increasing to $u_{t}$ as $\varepsilon
\rightarrow 0$, where $u_{t}$ is the solution of the BODE $%
u_{t}=x+\int_{t}^{T}\phi (u_{s})ds$.
\end{lemma}

\proof%
By comparison theorem \ref{compODE3} , we know that $u_{t}^{\varepsilon
_{1}}\geq u_{t}^{\varepsilon _{2}}$, for $\varepsilon _{1}\leq \varepsilon
_{2}$. So $u_{t}^{\varepsilon }\nearrow u_{t}$, for $t\in [0,T]$, as $%
\varepsilon \rightarrow 0$. Then the result follows easily from the
continuity of $\phi $ in $y$ and the boundedness of $u_{t}^{\varepsilon }$. $%
\square $

Now we prove a useful inequality.

\begin{lemma}
\label{ode-ineq}Let $u_{t}$ be the solution of BODE $u_{t}=x+\int_{t}^{T}%
\phi (u_{s})ds$, and $y_{t}$ which satisfies $y_{t}\geq x+\int_{t}^{T}\phi
(y_{s})ds$, on $[0,T]$. Then $u_{t}\leq y_{t}$, for $t\in [0,T]$.
\end{lemma}

\proof%
For any $\varepsilon >0$, $y_{T}\geq x>x-\varepsilon =u_{T}^{\varepsilon }$,
where $u_{t}^{\varepsilon }$ is the solution of (\ref{ode-epsilon}). Suppose
that there exists a $\overline{\tau }$, such that $u_{\overline{\tau }%
}^{\varepsilon }=y_{\overline{\tau }}$ and $y_{s}>u_{s}^{\varepsilon }$ on $[%
\overline{\tau },T]$. It follows from the monotonicity of $\phi $ on $y$
that
\[
y_{\overline{\tau }}\geq x+\int_{\overline{\tau }}^{T}\phi (y_{s})ds\geq
x+\int_{\overline{\tau }}^{T}\phi (u_{s}^{\varepsilon })ds>u_{\overline{\tau
}}^{\varepsilon },
\]
which is a contradiction. So $y_{t}>u_{t}^{\varepsilon }$ on $[0,T]$, for
any $\varepsilon >0$. Let $\varepsilon \rightarrow 0$, with lemma \ref
{ode-lim}, we have $u_{t}\leq y_{t}$, on $[0,T]$. $\square $

With the help of these Lemmas, we give the representation of the solution of
the reflected BODE.

\begin{proposition}
\label{rep}Under the assumptions A1, A2 and A3, assume that $%
(y_{t},k_{t})_{0\leq t\leq T}$ is a solution of the following reflected BODE
\begin{eqnarray}
y_{t} &=&x+\int_{t}^{T}\phi (y_{s})ds+k_{T}-k_{t},  \label{RODE2} \\
y_{t} &\geq &l_{t},\int_{0}^{T}(y_{s}-l_{s})dk_{s}=0.  \nonumber
\end{eqnarray}
Then we have for $t\in [0,T]$,
\[
y_{t}=\sup_{t\leq s\leq T}u_{t}^{s}=\sup_{t\leq s\leq T}[\int_{t}^{s}\phi
(y_{r})dr+x1_{\{s=T\}}+l_{s}1_{\{s<T\}}],
\]
where $(u_{t}^{s})_{0\leq t\leq s}$ is the solution of the BODE defined on $%
[0,s]$ with coefficient $\phi $ and terminal value $x1_{\{s=T\}}+l_{s}1_{%
\{s<T\}}$.
\end{proposition}

\proof%
For $t\in [0,T]$, since $y$ is a solution of the reflected BODE, with lemma
\ref{ode-ineq}, we get $y_{t}\geq u_{t}^{s}$, for $s\in [t,T]$. Denote
\[
D_{t}=\inf \{u\in [t,T],y_{u}=l_{u}\}\wedge T.
\]
Notice that $k$ is an increasing process and $%
\int_{0}^{T}(y_{s}-l_{s})dk_{s}=0$, then $k_{D_{t}}=k_{t}$. It follows that
\[
y_{t}=u_{t}^{D_{t}},
\]
which implies the first equality. For the second one, from (\ref{RODE2}), it
follows
\begin{eqnarray*}
y_{t} &=&y_{s}+\int_{t}^{s}\phi (y_{r})dr+k_{T}-k_{s} \\
&\geq &\int_{t}^{s}\phi (y_{r})dr+x1_{\{s=T\}}+l_{s}1_{\{s<T\}}.
\end{eqnarray*}
With the same $D_{t}$, we have
\[
y_{t}=\int_{t}^{D_{t}}\phi
(y_{r})dr+x1_{\{D_{t}=T\}}+l_{D_{t}}1_{\{D_{t}<T\}}.
\]
The proof is complete. $\square $

\begin{remark}
The function $H(p)=p(\alpha \gamma +\beta \ln p1_{\left[ 1,+\infty \right[
}(p))+\alpha \gamma 1_{\left( -\infty ,1\right) }(p)$ satisfies assumption
A3, the existence and uniqueness of the solution of non reflected BODE's
(see \cite{BH2005}). Consequently we have the result of the theorem \ref
{ode-main} relatively to $H$.
\end{remark}


\begin{thebibliography}{99}
\bibitem{BDHPS}  Briand, Ph., Delyon, B., Hu, Y., Pardoux, E. and Stoica, L.
(2003) $L_{p}$ solutions of BSDEs, \textit{Stochastic Process. Appl.} 108,
109-129.

\bibitem{BH2005}  Briand, Ph., Hu, Y. (2005) BSDE with quadratic growth and
unbounded terminal value, preprint.

\bibitem{E79}  N. El Karoui, (1979) Les aspects probabilistes du contr\^{o}%
le stochastique. Ecole d'\'{e}t\'{e} de Saint-Flour, \textit{Lecture Notes
in Math. }\textbf{876. (}Springer, Berlin), 73-238.

\bibitem{EKPPQ}  N.~El Karoui, C.~Kapoudjian, E.~Pardoux, S.~Peng and M.C.
Quenez, (1997) Reflected Solutions of Backward SDE and Related Obstacle
Problems for PDEs, \textit{Ann. Probab.} \textbf{25}, no 2, 702--737.

\bibitem{EPQ}  N.~El Karoui, S.~Peng and M.C. Quenez, (1997). Backward
stochastic differential equations in Finance. \textit{Math. Finance, }%
\textbf{7, }1-71\textbf{.}

\bibitem{K00}  M. Kobylanski, (2000) Backward stochastic differential
equations and partial differential equations with quadratic growth, \textit{%
Ann. Proba. }\textbf{28}, 558-602.\textbf{\ }

\bibitem{KLQT}  M. Kobylanski, J.P. Lepeltier, M.C. Quenez and S.Torres,
(2002) Reflected BSDE with superlinear quadratic coefficient. \textit{%
Probability and Mathematical Statistics. }Vol. \textbf{22}, 51-83.

\bibitem{LMX}  J.P. Lepeltier, A. Matoussi and M. Xu, (2005) Reflected BSDEs
under monotonicity and general increasing growth conditiond. \textit{%
Advanced in applied probability}, 2005, March, 1-26.

\bibitem{LS98}  J.P. Lepeltier and J. San Mart\'{i}n. (1998) Existence for
BSDE with superlinear-quadratic coefficient, \textit{Stochastic and
Stochastic reports,} \textbf{63}, 227-240.

\bibitem{LS04}  J.P. Lepeltier and J. San Mart\'{i}n, (2004) BSDE's with
continuous, monotonicity, and non-Lipschitz in $z$ coefficient. Preprint

\bibitem{M97}  A. Matoussi, (1997) Reflected solutions of backward
stochastic differential equations with continuous coefficient, \textit{%
Statistic \& Probality Letters} \textbf{34}, 347-354.

\bibitem{PP90}  E. Pardoux and S. Peng, (1990) Adapted solutions of Backward
Stochastic Differential Equations. \textit{Systems Control Lett. }\textbf{14,%
} 51-61.

\bibitem{X2004}  Xu, M. (2004)Reflected BSDE with monotonicity and general
increasing in $y$, and non-Lipschitz conditions in $z$. Submitted to \textit{%
Stochastic Process. Appl.}.

\bibitem{RY}  D. Revuz and M. Yor, (1991) Continuous martingales and
Brownian motion (Springer, Berlin).
\end{thebibliography}
\end{document}